\documentclass[12pt]{article}
\usepackage[T2A]{fontenc}
\usepackage[cp1251]{inputenc}
\usepackage[english,russian]{babel}
\usepackage{amsmath}
\usepackage{amsfonts,amssymb}
\usepackage{graphicx}
\newtheorem{theorem}{Теорема}
\newtheorem{theorem2}{Theorem}
\begin{document}

\title{On Evaluation of Riesz Constants for Systems of Shifted Gaussians}

\author{E.\,A. Kiselev,\\ Voronezh State University, Voronezh, Russia,\\ evg-kisel2006@yandex.ru;\\
L.\,A. Minin,\\ Voronezh State University, Voronezh, Russia,\\ mininla@mail.ru;\\
I.\,Ya. Novikov,\\ Voronezh State University,Voronezh, Russia,\\ igor.nvkv@gmail.com;\\
S.\,M. Sitnik,\\ Voronezh Institute of the Ministry of Internal Affairs, Voronezh, Russia,\\ pochtasms@gmail.com.}

\date{}
\maketitle

\newpage
\begin{abstract}

\begin{center}
ANNOTATION\\
\end{center}

In this paper we study single--parametric systems of integer shifts of Gauss and Lorenz functions. In case of Lorenz system we explicitly calculate nod functions and prove that it tends to sinc function in limit. For both Gauss and  Lorenz systems and corresponding nod functions we explicitly calculate Riesz constants via trigonometric, hyperbolic and Jacobi theta--functions, also limit behavior of this values is found depending on parameters.

A special result is a sharp monotonicity property proved for a special ratio of  Jacobi theta--functions which is important in different applications.

At the end of the paper some open problems and further applications are considered.
\end{abstract}

\newpage

\section{Key definitions and main results}

\

\textsc{Definition 1.} (\cite{KS}, \cite{NPS}).
A set of functions $\varphi_k(t) \in L_{2}(\mathbf{R}), \, k \in \mathbf{Z},$ is the Riesz system if a pair of positive constants $A,B$ exists such that for every sequence $c\in l_{2}$ the next estimates are valid

\begin{equation*}
A\|c\|^{2}_{l_2}\leq\left\|\sum_{k=-\infty}^\infty c_{k} \varphi_{k}(t)\right\|^{2}_{L_{2}}\leq B\|c\|^{2}_{l_2} \, .
\label{KR1}
\end{equation*}

The max $A$ is called the lower Riesz constant and min $B$ the upper Riesz constant. The ratio of Riesz constants $\frac{B}{A}$ is an important system characteristic especially for numerical evaluations.

Many approximation methods are based on multiplicative shifts of a single function $\varphi$:

\begin{equation*}
\sum_k f_k \varphi(a_k t),
\end{equation*}
or additive shifts of  $\varphi$:
\begin{equation*}
\sum_k f_k \varphi(t-b_k),
\end{equation*}
(Fourier and Kaptein series, Wittaker--Shannon--Kotelnikov, Streng--Fix, Rvachev brothers, wavelets, frames and so on, cf. (\cite{{KS},{NPS},{Chui},{Ch}}))

\textsc{Definition 2.} The function $\widetilde{\varphi}(t)$ defined as
\begin{equation}
\widetilde{\varphi}(t)=\sum_{k=-\infty}^\infty d_{k}\varphi(t-k),
\label{Int3}
\end{equation}
is called a basic nod function if the next is valid
\begin{equation}
\widetilde{\varphi}(m)=\delta_{0m}, \, m\in\mathbb{Z},
\label{Int4}
\end{equation}
$\delta_{0m}$ -- Kronecker symbol.
This function solves interpolation problem because for the function
\begin{equation*}
\widetilde{f}(t)=\sum_{n=-\infty}^\infty f(n) \widetilde{\varphi}(t-n)
\end{equation*}
it is valid for integer values $\widetilde{f}(m)=f(m)$ при $m\in\mathbb{Z}$.

In this paper we study integer additive shifts of Gauss $\varphi(t)=\exp\left(-\frac{t^{2}}{2\sigma^{2}}\right)$ and Cauchy--Lorentz functions
$\varphi(t)=\frac{\sigma^{2}}{\sigma^{2}+t^{2}}$.

We use index $G$ for Gauss and $L$ for Cauchy--Lorentz connected values depending on a parameter $\sigma>0$..

\begin{table}[h]
\caption{Table 1: Basis notations}
\label{tabular:mainnote}
\begin{center}
\begin{tabular}{|c|c|c|}  \hline
Given function & Gauss function & Cauchy--Lorentz function  \\
$\varphi(t)$ & $\varphi_{G}(t,\sigma)=\exp{\left(-\frac{t^{2}}{2\sigma^{2}}\right)}$ & $\varphi_{L}(t,\sigma)=\frac{\sigma^{2}}{\sigma^{2}+t^{2}}$  \\  \hline
Riesz constants for  &  &  \\
shifted system $\varphi(t-k)$  & $A_{G}(\sigma), B_{G}(\sigma)$ & $A_{L}(\sigma), B_{L}(\sigma)$  \\  \hline
Series  coefficients $d_k$ & & \\
for basic nod function & $d_{G,k}(\sigma)$ & $d_{L,k}(\sigma)$  \\  \hline
Basis nod function $\widetilde{\varphi}(t)$ & $\widetilde{\varphi}_{G}(t,\sigma)$ & $\widetilde{\varphi}_{L}(t,\sigma)$ \\  \hline
Riesz constants for &  &  \\
shifts of nod function  & $\widetilde{A}_{G}(\sigma), \widetilde{B}_{G}(\sigma)$ & $\widetilde{A}_{L}(\sigma), \widetilde{B}_{L}(\sigma)$  \\  \hline
\end{tabular}
\end{center}
\end{table}

It was proved in \cite {ZhKMS} that integer shifts of the Gauss function form the Riesz system with explicit constants

\begin{equation*}
A_{G}(\sigma)=\sigma\sqrt{\pi}\vartheta_{3}\left(\frac{\pi}{2},q\right), \,
B_{G}(\sigma)=\sigma\sqrt{\pi}\vartheta_{3}\left(0,q \right), \,
q=\exp\left({-\frac{1}{4\sigma^2}}\right),
\end{equation*}
and $\vartheta_3(t,q)$ -- the third Jacobi theta--function \cite{UW}
\begin{equation}
\vartheta_3(t,q)=\sum_{k=-\infty}^\infty q^{k^{2}}e^{2ikt},|q|<1.
\label{Teta3}
\end{equation}

In this paper we prove also explicit formulas for Cauchy--Lorentz system.

\begin{theorem2}
Integer shifts of Cauchy--Lorentz function form the Riesz system with explicit Riesz constants
\begin{equation*}
A_L(\sigma)=\frac{\sigma^{2}\pi^{2}}{\sh(2\sigma\pi)}, \;
B_L(\sigma)=\frac{\sigma^{2}\pi^{2} \ch(2\sigma\pi)}{\sh(2\sigma\pi)}.
\end{equation*}
\end{theorem2}
\newpage

This is a table of Riesz constants for both cases.

\begin{table}[h]
\caption{Table 2: Riesz constants for shifted systems}
\centerline{of Gauss and Cauchy--Lorentz function}
\label{tabular:znachRiss1}
\begin{center}
\begin{tabular}{|c|c|c|c|c|c|c|}  \hline
$\sigma$ & $A_{G}(\sigma)$ & $B_{G}(\sigma)$ & $B_{G}(\sigma)/A_{G}(\sigma)$ & $A_{L}(\sigma)$ & $B_{L}(\sigma)$ & $B_{L}(\sigma)/A_{L}(\sigma)$ \\  \hline
0{.}2 & $0.353$ & $0.356$ & $1.01$
 & $0.245$ & $0.464$ & $1.90$ \\  \hline
0{.}4 & $0.415$ & $1.009$ & $2.43$
 & $0.258$ & $1.600$ & $6.21$ \\  \hline
0{.}6 & $0.130$ & $2.262$ & $17.46$
 & $0.164$ & $3.557$ & $21.70$ \\  \hline
1{.}0 & $6.45\cdot10^{-4}$ & $6.283$ & $9.67\cdot10^{3}$
 & $0.037$ & $9.870$ & $267.75$ \\  \hline
2{.}0 & $3.60\cdot10^{-16}$ & $25.13$ & $6.98\cdot10^{16}$
 & $2.75\cdot10^{-4}$ & $39.48$ & $1.43\cdot10^{5}$ \\  \hline
3{.}0 & $3.00\cdot10^{-37}$ & $56.55$ & $1.88\cdot10^{38}$
 & $1.16\cdot10^{-6}$ & $88.83$ & $7.68\cdot10^{7}$ \\  \hline
4{.}0 & $5.28\cdot10^{-67}$ & $100.53$ & $1.91\cdot10^{68}$
 & $3.84\cdot10^{-9}$ & $157.91$ & $4.11\cdot10^{10}$ \\  \hline
5{.}0 & $2.18\cdot10^{-105}$ & $157.08$ & $7.19\cdot10^{106}$
 & $1.12\cdot10^{-11}$ & $246.74$ & $2.20\cdot10^{13}$ \\  \hline
\end{tabular}
\end{center}
\end{table}

As for basic nod function for Gauss system its explicit value was found in \cite[p. 152]{Maz}. Coefficient formulas for $d_k$ in (\ref {Int3}) are:

\begin{equation*}
d_{G,k}(\sigma) =\frac{1}{C\left(\sigma \right)} \cdot \exp \left(\frac{k^{2} }{2\sigma ^{2} } \right)\cdot \sum \limits _{r=\left|k\right|}^{\infty }\left(-1\right)^{r} \cdot \exp \left(-\frac{\left(r+0.5\right)^{2} }{2\sigma ^{2}} \right),
\end{equation*}
где
\begin{equation*}
C\left(\sigma \right)=
\sum \limits _{r=-\infty }^{\infty }\left(4r+1\right)\cdot
\exp \left(-\frac{\left(2r+0.5\right)^{2}}{2\sigma ^{2}} \right).
\end{equation*}

We prove similar result for Cauchy--Lorentz system.

\begin{theorem2} The next formula for coefficients of basic nod function in case of Cauchy--Lorentz system is valid:
\begin{equation*}
d_{L,k}(\sigma)=\frac{(-1)^{k}\sh(\sigma\pi)}{\sigma\pi^{2}}\int \limits_{0}^{\pi}\frac{\cos(kt)}{\ch(\sigma t)}dt.
\end{equation*}
\end{theorem2}

Riesz constants are different for systems of shifts of nod functions and generating functions though their subspaces are the same.

\begin{theorem2}
The next limit formulas are valid

\begin{equation*}
\lim_{\sigma\rightarrow\infty}\widetilde{A}_{G}(\sigma)=
\lim_{\sigma\rightarrow\infty}\widetilde{A}_{L}(\sigma)=
\frac{1}{2},
\end{equation*}
\begin{equation*}
\lim_{\sigma\rightarrow\infty}\widetilde{B}_{G}(\sigma)=
\lim_{\sigma\rightarrow\infty}\widetilde{B}_{L}(\sigma)=1.
\end{equation*}
\end{theorem2}

In \cite {Zhur} for Gauss function the estimate for $\widetilde{B}_{G}(\sigma)$ was proved but for
$\widetilde{A}_{G}(\sigma)$ just a weaker result was proved
$$\overline{\lim_{\sigma\rightarrow\infty}}\widetilde{A}_{G}(\sigma)\leq\frac{1}{2}.
$$

Now consider basic nod functions for $\sigma \rightarrow \infty$.
It depends on $sinc$--function, cf. \cite{Ste}. It is known \cite{Sivakumar} that for Gauss system
\begin{equation*}
\widetilde{\varphi}_{G}(t,\sigma) \xrightarrow[\sigma \rightarrow \infty]{L_2(\mathbb{R})} \mathrm{sinc}(\pi t),
\end{equation*}
\begin{equation*}
\mathrm{sinc}\left(\pi t\right)=\frac{\sin \pi t}{\pi t}.
\end{equation*}
We prove similar result for Cauchy--Lorentz system.

\begin{theorem2}
The next limit formulas are valid
\begin{equation*}
\widetilde{\varphi}_{L}(t,\sigma) \xrightarrow[\sigma \rightarrow \infty]{L_2(\mathbb{R})} \mathrm{sinc}\left(\pi t\right).
\end{equation*}
\end{theorem2}

To prove limit results for the Gauss system the next important monotonicity property for Jacobi $\vartheta$--function is proved.

\begin{theorem2}
Define the function $P(t)$  as
\begin{equation}
P(t)=\frac{\vartheta_3 (t, q)}
     {\left[ \vartheta_3 \left( t,q^2 \right) \right]^2}, \
      t \in [0, \pi]; \, q=\exp \left( -\frac{1}{4\sigma^2} \right).
\label{PGT1}
\end{equation}
Then $P(t)$ is monotonic decreasing for $t\in(0, \frac {\pi} {2})$ and monotonic increasing for $(t\in\frac {\pi} {2}, \pi)$.
\end{theorem2}

Inequalities and monotonicity properties for ratios of $\vartheta$--functions are important in different fields and has the long history. They are the key estimates for A.A.~Gonchar problem in complex theory and V.N.~Dubinin method  \cite{Sol} and other applications, cf. \cite{Dixit} -- \cite{Schief}. On inequalities
for similiar special functions cf. \cite{S1}--\cite{S6}.
\section{Notes to applications}

\

For shifted systems and approximations there are too many references, we do not try to include even a small part of it and restrict ourselves to only small number of ones we used directly.

In papers \cite{Vin1}--\cite{Vin2}  more general expansions are considered by the functions
$$
\sum\limits_{k=-\infty}^{\infty} f_k \frac{a}{b+(x-k)^c},
$$
with some positive constants $a,b,c$. In this paper we consider a case $c=2$. 

\bigskip

Many applications are in \cite{Maz},\cite{Ch},\cite{BS} and in bibliographies of these references.

\bigskip

Above it let us mention some more topics.

\bigskip

Problem 1. Prove that finite dimensional solutions studied in  \cite{ST1}--\cite{ST5} for nod function  (\ref{Int3}) converge to the exact solution from \cite{Maz} when linear system dimension tends to infinity.

\bigskip

Problem 2. Prove symmetry and sign changing for coefficients in (\ref{Int3}), some results are proved in \cite{MSU}, \cite{KMNS1}--\cite{KMNS2}.

\bigskip

Problem 3. Study expansions with half--integer shifts.

\bigskip

Problem 4. Study connections with theory of coherent states  \cite{Per}. It is connected with famous von Neumann hypothesis which is not clear completely still, cf. \cite{SMZ}.

\bigskip

Problem 5. Study  numeric methods based on Discreet Fourier Transform \cite{SMZ},\cite{S17}.

\bigskip

Problem 6. Gabor frames, tight frames \cite{Leb}, \cite{Chr2}, \cite{Chr3}.

\bigskip

Problem 7. Further use of inequalities for Jacobi $\vartheta$ and other special functions, cf. \cite{MS1}, \cite{S1}-\cite{S6}.

\bigskip

Problem 8. Applications and connections with transmutation theory, cf.   \cite{S7}-\cite{S16}.

\bigskip

\bigskip

The russian version is ended with a list of some unsolved problems on page 27.

\newpage

\section{Русский вариант статьи:}
\begin{center}
\Large
Е.А. Киселёв, Л.А. Минин, И.Я. Новиков \\
(Воронежский государственный университет),\\
С.М. Ситник \\
(Воронежский институт МВД России).\\
\bigskip
\Large\bf{
Вычисление констант Рисса\\
для системы целочисленных сдвигов функции Гаусса.}
\bigskip
--------------------------------------------------------
\bigskip

\large
Краткое содержание работы.
\end{center}

\

В работе изучаются однопараметрические семейства целочисленных сдвигов функций Гаусса и Лоренца.
В случае функции Лоренца построены в явном виде узловые функции и показано, что предельным значением для них является функция отсчетов. Для систем сдвигов, порожденных как функциями Гаусса и Лоренца, так и связанными с ними узловыми функциями, получены явные выражения для констант Рисса и изучено поведение этих констант в зависимости от параметра.

Установленная при доказательстве результатов работы монотонность  одного специального отношения тета--функций Якоби имеет самостоятельное значение.

\newpage

\section{Введение}

\

Время до 1950-х годов можно назвать господством ортогональных систем функций в математике, которое длилось почти два с половиной века \cite{KS}. После этого рубежа началось время широкого использования  неортогональных систем, поскольку они лучше приспособлены для описания локальных свойств изучаемых функций, а не их упрощенных средних характеристик. Обычно платой за неортогональность является неустойчивость численных методов решения поставленных задач. Мерой неустойчивости служит соотношение констант Рисса.

\textsc{Определение 1.} (\cite{KS}, \cite{NPS}).
Функции $\varphi_k(t) \in L_{2}(\mathbf{R}), \, k \in \mathbf{Z},$ образуют систему Рисса, если существуют
положительные константы $A$ и $B$ такие, что для любой последовательности коэффициентов
$c\in l_{2}$ выполнена двусторонняя оценка
\begin{equation*}
A\|c\|^{2}_{l_2}\leq\left\|\sum_{k=-\infty}^\infty c_{k} \varphi_{k}(t)\right\|^{2}_{L_{2}}\leq B\|c\|^{2}_{l_2} \, .
\label{KR1}
\end{equation*}
Наибольшая из величин $A$ называется нижней константой Рисса, наименьшая  из величин $B$  -- верхней константой Рисса. Для ортонормированных систем функций обе константы равны 1.
В случае конечного набора линейно независимых функций $\varphi_{k}(t)$ аналогами констант Рисса являются минимальное и максимальное собственные значения матрицы Грама, образованной попарными скалярными произведениями этих функций \cite{KS}. Число обусловленности матрицы Грама равно отношению максимального собственного значения к минимальному. Если оно велико, то матрица называется плохо обусловленной и при работе с ней требуется применять специальные приемы для обеспечения устойчивости вычислений \cite{BZhK}.

Широкий класс разложений в функциональные ряды может быть охарактеризован как аппроксимации по мультипликативным или аддитивным сдвигам одной функции $\varphi (t)$.
Наиболее известными представителями семейства разложений по мультипликативным сдвигам вида
\begin{equation*}
\sum_k f_k \varphi(a_k t),
\end{equation*}
где $f_k$ -- искомые коэффициенты разложения заданной функции $f(t)$,\\ $a_k$ -- заданный набор сдвигов (масштабирующих коэффициентов), являются ряды Фурье, ряды Каптейна и Шлемильха по функциям Бесселя \cite{JP}.

Второе семейство -- это аппроксимации аддитивными сдвигами вида
\begin{equation*}
\sum_k f_k \varphi(t-b_k),
\end{equation*}
где $b_k$ -- заданный набор сдвигов. Сюда относятся разложения с использованием формулы Уиттекера--Шеннона--Котельникова, сплайнов Стренга--Фикса,  функций Рвачёвых и многие другие.

На комбинациях двух указанных семейств разложений построены теории всплесков \cite{{KS},{NPS},{Chui}} и фреймов \cite{{Ch}}.

\textsc{Определение 2.} Функция $\widetilde{\varphi}(t)$, являющаяся линейной комбинацией
системы целочисленных аддитивных сдвигов $\varphi_{k}(t)=\varphi(t-k), \, k \in \mathbf{Z}$,
\begin{equation}
\widetilde{\varphi}(t)=\sum_{k=-\infty}^\infty d_{k}\varphi(t-k),
\label{Int3}
\end{equation}
называется узловой функцией, если для нее выполнена система равенств
\begin{equation}
\widetilde{\varphi}(m)=\delta_{0m}, \, m\in\mathbb{Z},
\label{Int4}
\end{equation}
где $\delta_{0m}$ -- символ Кронекера.

С помощью узловой функции решается задача интерполяции на бесконечной равномерной сетке, поскольку функция
\begin{equation*}
\widetilde{f}(t)=\sum_{n=-\infty}^\infty f(n) \widetilde{\varphi}(t-n)
\end{equation*}
совпадает с функцией $f(t)$ в целых узлах, т.е. $\widetilde{f}(m)=f(m)$ при $m\in\mathbb{Z}$.

В данной статье рассматриваются системы целочисленных аддитивных сдвигов функции Гаусса
$\varphi(t)=\exp\left(-\frac{t^{2}}{2\sigma^{2}}\right)$ и функции Лоренца
 $\varphi(t)=\frac{\sigma^{2}}{\sigma^{2}+t^{2}}$.
Для исходных систем, а также для систем, порожденных сдвигами соответствующих узловых функций, получены явные аналитические формулы для констант Рисса и изучено их предельное поведение при $\sigma\rightarrow\infty$.

Название  функция Лоренца мы используем для определённости. Более точно, эта функция возникла у Коши как плотность соответствующего вероятностного распределения, в физических приложения используются названия функция Лоренца или функция Брейта--Вигнера.

Использование квадратичных экспонент --- или функций Гаусса --- имеет давнюю историю. Впервые они возникли у Френеля, который предложил использовать квадратичные экспоненты вместо обычных плоских волн, так как только на этом пути получалось совпадающее с опытом объяснение явления  Фраунгофера; при этом математика обогатилась известными интегралами Френеля. В середине 20--го века именно разложения по целочисленным сдвигам квадратичных экспонент стали тем математическим аппаратом, при помощи которого Габор теоретически обосновал, а затем практически реализовал идеи голографии. В начале 21--го века компьютерный пакет "Gaussian"\,, основанный на использовании разложений по квадратичным экспонентам, был оценён Нобелевской премией, что пока является единственным случаем присуждения этой премии за компьютерную программу; причиной явилось решение с помощью данного пакета давно стоявших задач на стыке физики, химии и биологии. Кроме приведённых примеров дискретных разложений по квадратичным экспонентам используются и непрерывные разложения. Эти разложения имеют многочисленные названия (дробное преобразование Фурье, квадратичное преобразование Фурье, преобразование Габора, преобразование Фурье--Френеля, преобразование Вейерштрасса, полугруппа Орстейна--Уленбека, преобразование Винера, преобразование Бора--Френеля, интегральное преобразование Гаусса и т.д.), см., например, \cite{FFS}. В теории чисел ряды и суммы квадратичных экспонент --- суммы Гаусса --- являются одним из наиболее важных и изучаемых объектов \cite{SG1}--\cite{SG2}.

Использование функций Коши--Лоренца первоначально определяется тем, что при выборе соответствующих параметров они практически не отличаются на графиках от функций Гаусса, но имеют несколько более простую структуру.
Для приближения функций Гаусса с конкретными параметрами функциями Коши--Лоренца разумно выбирать последние как аппроксимации Паде порядка [0/2] для квадратичных экспонент.

\section{Основные результаты}

\

Для функции Гаусса в дальнейшем используется индекс $G$, для функции Лоренца -- индекс $L$. Все изучаемые в работе  величины рассматриваются как функции, зависящие от параметра $\sigma>0$. В таблице
\ref{tabular:mainnote} представлены введенные нами обозначения.

\begin{table}[h]
\caption{Основные обозначения}
\label{tabular:mainnote}
\begin{center}
\begin{tabular}{|c|c|c|}  \hline
Исходная функция & Функция Гаусса & Функция Лоренца  \\
$\varphi(t)$ & $\varphi_{G}(t,\sigma)=\exp{\left(-\frac{t^{2}}{2\sigma^{2}}\right)}$ & $\varphi_{L}(t,\sigma)=\frac{\sigma^{2}}{\sigma^{2}+t^{2}}$  \\  \hline
Константы Рисса для  &  &  \\
системы сдвигов $\varphi(t-k)$  & $A_{G}(\sigma), B_{G}(\sigma)$ & $A_{L}(\sigma), B_{L}(\sigma)$  \\  \hline
Коэффициенты $d_k$ ряда & & \\
для узловой функции & $d_{G,k}(\sigma)$ & $d_{L,k}(\sigma)$  \\  \hline
Узловая функция $\widetilde{\varphi}(t)$ & $\widetilde{\varphi}_{G}(t,\sigma)$ & $\widetilde{\varphi}_{L}(t,\sigma)$ \\  \hline
Константы Рисса для &  &  \\
сдвигов узловой функции  & $\widetilde{A}_{G}(\sigma), \widetilde{B}_{G}(\sigma)$ & $\widetilde{A}_{L}(\sigma), \widetilde{B}_{L}(\sigma)$  \\  \hline
\end{tabular}
\end{center}
\end{table}

В статье \cite {ZhKMS} показано, что целочисленные сдвиги функции Гаусса образуют систему Рисса с константами
\begin{equation*}
A_{G}(\sigma)=\sigma\sqrt{\pi}\vartheta_{3}\left(\frac{\pi}{2},q\right), \,
B_{G}(\sigma)=\sigma\sqrt{\pi}\vartheta_{3}\left(0,q \right), \,
q=\exp\left({-\frac{1}{4\sigma^2}}\right),
\end{equation*}
где $\vartheta_3(t,q)$ -- третья тета-функция Якоби  \cite{UW}
\begin{equation}
\vartheta_3(t,q)=\sum_{k=-\infty}^\infty q^{k^{2}}e^{2ikt},|q|<1.
\label{Teta3}
\end{equation}

В настоящей работе получены также явные выражения для констант Рисса в случае функции Лоренца.
\begin{theorem}
Целочисленные сдвиги функции Лоренца образуют систему Рисса с константами
\begin{equation*}
A_L(\sigma)=\frac{\sigma^{2}\pi^{2}}{\sh(2\sigma\pi)}, \;
B_L(\sigma)=\frac{\sigma^{2}\pi^{2} \ch(2\sigma\pi)}{\sh(2\sigma\pi)}.
\end{equation*}
\end{theorem}

Численные значения констант Рисса для обоих случаев представлены таблице \ref{tabular:znachRiss1}.
Все значащие цифры -- верные (с точностью до округления).

\begin{table}[h]
\caption{Значения констант Рисса для систем сдвигов,}
\centerline{порожденных функцией Гаусса и функцией Лоренца}
\label{tabular:znachRiss1}
\begin{center}
\begin{tabular}{|c|c|c|c|c|c|c|}  \hline
$\sigma$ & $A_{G}(\sigma)$ & $B_{G}(\sigma)$ & $B_{G}(\sigma)/A_{G}(\sigma)$ & $A_{L}(\sigma)$ & $B_{L}(\sigma)$ & $B_{L}(\sigma)/A_{L}(\sigma)$ \\  \hline
0{.}2 & $0.353$ & $0.356$ & $1.01$
 & $0.245$ & $0.464$ & $1.90$ \\  \hline
0{.}4 & $0.415$ & $1.009$ & $2.43$
 & $0.258$ & $1.600$ & $6.21$ \\  \hline
0{.}6 & $0.130$ & $2.262$ & $17.46$
 & $0.164$ & $3.557$ & $21.70$ \\  \hline
1{.}0 & $6.45\cdot10^{-4}$ & $6.283$ & $9.67\cdot10^{3}$
 & $0.037$ & $9.870$ & $267.75$ \\  \hline
2{.}0 & $3.60\cdot10^{-16}$ & $25.13$ & $6.98\cdot10^{16}$
 & $2.75\cdot10^{-4}$ & $39.48$ & $1.43\cdot10^{5}$ \\  \hline
3{.}0 & $3.00\cdot10^{-37}$ & $56.55$ & $1.88\cdot10^{38}$
 & $1.16\cdot10^{-6}$ & $88.83$ & $7.68\cdot10^{7}$ \\  \hline
4{.}0 & $5.28\cdot10^{-67}$ & $100.53$ & $1.91\cdot10^{68}$
 & $3.84\cdot10^{-9}$ & $157.91$ & $4.11\cdot10^{10}$ \\  \hline
5{.}0 & $2.18\cdot10^{-105}$ & $157.08$ & $7.19\cdot10^{106}$
 & $1.12\cdot10^{-11}$ & $246.74$ & $2.20\cdot10^{13}$ \\  \hline
\end{tabular}
\end{center}
\end{table}

Перейдем к узловым функциям.
В случае функции Гаусса аналитическое выражение для коэффициентов $d_k$ из формулы (\ref {Int3}) приведено в
\cite[с. 152]{Maz}:
\begin{equation*}
d_{G,k}(\sigma) =\frac{1}{C\left(\sigma \right)} \cdot \exp \left(\frac{k^{2} }{2\sigma ^{2} } \right)\cdot \sum \limits _{r=\left|k\right|}^{\infty }\left(-1\right)^{r} \cdot \exp \left(-\frac{\left(r+0.5\right)^{2} }{2\sigma ^{2}} \right),
\end{equation*}
где
\begin{equation*}
C\left(\sigma \right)=
\sum \limits _{r=-\infty }^{\infty }\left(4r+1\right)\cdot
\exp \left(-\frac{\left(2r+0.5\right)^{2}}{2\sigma ^{2}} \right).
\end{equation*}

Аналогичное утверждение нами доказано и для случая функции Лоренца.
\begin{theorem} Для коэффициентов узловой функции, построенной по системе сдвигов функции Лоренца, справедлива формула
\begin{equation*}
d_{L,k}(\sigma)=\frac{(-1)^{k}\sh(\sigma\pi)}{\sigma\pi^{2}}\int \limits_{0}^{\pi}\frac{\cos(kt)}{\ch(\sigma t)}dt.
\end{equation*}
\end{theorem}

Подпространства, порожденные сдвигами узловой функции, совпадают с подпространствами, порожденными сдвигами исходной функции. Но константы Рисса различаются, поскольку базисы разные.

\begin{theorem}
Справедливы следующие предельные соотношения:

\begin{equation*}
\lim_{\sigma\rightarrow\infty}\widetilde{A}_{G}(\sigma)=
\lim_{\sigma\rightarrow\infty}\widetilde{A}_{L}(\sigma)=
\frac{1}{2},
\end{equation*}
\begin{equation*}
\lim_{\sigma\rightarrow\infty}\widetilde{B}_{G}(\sigma)=
\lim_{\sigma\rightarrow\infty}\widetilde{B}_{L}(\sigma)=1.
\end{equation*}

\end{theorem}
Соответственно, предел отношения верхней константы Рисса к нижней для узловых функций равен 2.

Случай функции Гаусса рассмотрен в \cite {Zhur}. Утверждение для $\widetilde{B}_{G}(\sigma)$ там доказано, а для  $\widetilde{A}_{G}(\sigma)$ получен более слабый результат:
$$\overline{\lim_{\sigma\rightarrow\infty}}\widetilde{A}_{G}(\sigma)\leq\frac{1}{2}.
$$

Рассмотрим теперь поведение узловых функций при $\sigma \rightarrow \infty$.
Это предельное поведение определяется известной функцией отсчётов, или $sinc$--функцией. Теория разложений по $sinc$--функциям, включающая знаменитую формулу Уиттекера--Котельникова--Шеннона, подробно разработана, см., например, \cite{Ste}.
Известно \cite{Sivakumar}, что в случае функции Гаусса
\begin{equation*}
\widetilde{\varphi}_{G}(t,\sigma) \xrightarrow[\sigma \rightarrow \infty]{L_2(\mathbb{R})} \mathrm{sinc}(\pi t),
\end{equation*}
где функция отсчетов определяется так:
\begin{equation*}
\mathrm{sinc}\left(\pi t\right)=\frac{\sin \pi t}{\pi t}.
\end{equation*}
Нами установлено аналогичное утверждение и для функции Лоренца.
\begin{theorem}
Справедливо предельное соотношение
\begin{equation*}
\widetilde{\varphi}_{L}(t,\sigma) \xrightarrow[\sigma \rightarrow \infty]{L_2(\mathbb{R})} \mathrm{sinc}\left(\pi t\right).
\end{equation*}
\end{theorem}
Функция отсчетов порождает ортонормированную систему целочисленных сдвигов. Следовательно, предел нижней константы Рисса для обеих изучаемых узловых функций не равен нижней константе от предельной функции.

\section{Доказательства}

\

Основные инструменты при изучении систем сдвигов -- это преобразование Фурье
\begin{equation}
\widehat{\varphi}(\omega)=\frac{1}{\sqrt{2\pi}}\int \limits_{-\infty}^{\infty}\varphi(t)e^{-i\omega t}dt,
\label{FT}
\end{equation}
являющееся при таком выборе нормировочного множителя унитарным оператором в $L_2(\mathbf{R})$,
и формула суммирования Пуассона \cite{UW}
\begin{equation}
\sum \limits_{k=-\infty}^{\infty} \varphi(k)\cdot e^{-ik\omega}=\sqrt{2\pi}\sum \limits_{k=-\infty}^{\infty} \widehat{\varphi}(\omega+2 \pi k). \label{puasson}
\end{equation}
Выпишем образы Фурье по переменной $t$ для функций Гаусса и Лоренца:
\begin{equation}
\widehat{\varphi}_{G}(\omega,\sigma)=\sigma \exp{\left(-\frac{\sigma^{2}\omega^{2}}{2}\right)}, \;
\widehat{\varphi}_{L}(\omega,\sigma)=\sigma\sqrt{\frac{\pi}{2}}e^{-\sigma|\omega|}.
\label{FGauss_FLorentz}
\end{equation}
Формула (\ref{puasson}) для двух рассматриваемых функций имеет вид:
\begin{equation}
\sum \limits_{k=-\infty}^{\infty} \exp{\left(-\frac{k^{2}}{2\sigma^{2}}\right)}\cdot e^{-ik\omega}=\sigma\sqrt{2\pi}\sum \limits_{k=-\infty}^{\infty} \exp{\left(-\frac{\sigma^{2}(\omega+2\pi k)^{2}}{2}\right)},
\label{PGauss}
\end{equation}
\begin{equation}
\sum \limits_{k=-\infty}^{\infty} \frac{\sigma^{2}}{\sigma^{2}+k^{2}} \cdot e^{-ik\omega}=\sigma\pi\sum \limits_{k=-\infty}^{\infty} e^{-\sigma|\omega+2\pi k|}.
\label{PLorentz}
\end{equation}

Из равенств (\ref{Int3}) и (\ref{Int4}) получается следующая бесконечная система уравнений относительно неизвестных $d_{k}$:
\begin{equation}
\sum \limits_{k=-\infty}^{\infty}d_{k}\varphi(m-k)=\delta_{0m},\label{Intd}
\end{equation}
которая представляет собой систему типа свертки. Рассмотрим ряды Фурье с коэффициентами $d_{k}$ и $\varphi(k)$
\begin{equation*}
D(t)=\sum_{k=-\infty}^\infty d_{k}e^{-ikt}, \; \Phi(t)=\sum_{k=-\infty}^\infty \varphi(k)e^{-ikt}.
\end{equation*}
Система уравнений (\ref {Intd}) может быть записана в виде функционального равенства  \cite[с. 151]{Maz}
\begin{equation}
D(t) \cdot \Phi(t)=1.\label{Int8}
\end{equation}
Следовательно, для нахождения коэффициентов $d_{k}$ нужно разложить в ряд Фурье функцию $1/\Phi(t)$.

Ряд Фурье $D(t)$ называется маской \cite{NPS} (или символом \cite{Chui}) последовательности
$\{d_{k} \}_{k\in\mathbb{Z}}$.
В случае функции Гаусса маски $D(t)$ и $\Phi(t)$ будем обозначать как $D_{G}(t,\sigma)$ и $\Phi_{G}(t,\sigma)$,
в случае функции Лоренца --  как $D_{L}(t,\sigma)$ и $\Phi_{L}(t,\sigma)$.

\subsection{Доказательство теоремы 1}
\

Для систем целочисленных сдвигов одной заданной функции известно следующее утверждение.
\begin{theorem} \textup{(см. \cite{NPS}.)}
\label{theorem:NRS}
Пусть $\varphi\in L_{2}(\mathbb{R})$. Для того чтобы система функций $\varphi(t-k), k\in \mathbb{Z}$, являлась системой Рисса с постоянными $A, B$, необходимо и достаточно, чтобы для почти всех $\omega\in \mathbb{R}$ выполнялось соотношение
\begin{equation}
A\leq 2\pi \sum_{k\in\mathbb{Z}}\left|\widehat{\varphi}(\omega+2\pi k)\right|^2\leq B .\label{KRGL0}
\end{equation}
\end{theorem}

Сумма в неравенстве (\ref{KRGL0}) представляет собой периодическую функцию от $\omega$ с периодом $2\pi$.
Поэтому, если обозначить
\begin{equation}
P(\omega)=2\pi\sum_{k\in\mathbb{Z}}\left|\widehat{\varphi}(\omega+2\pi k)\right|^2, \label{PKR0}
\end{equation}
то константы Рисса находятся с помощью следующих соотношений
\begin{equation*}
A=\inf \limits_{0\leq\omega\leq2\pi} P(\omega), \,
B=\sup \limits_{0\leq\omega\leq2\pi} P(\omega).
\end{equation*}
В случае функции Гаусса $P(\omega)$ будем обозначать как $P_{G}(\omega,\sigma)$, а
соответствующий ряд для узловой функции -- как $\widetilde{P}_{G}(\omega,\sigma)$. Аналогично определяются для функции Лоренца величины  $P_{L}(\omega,\sigma)$ и  $\widetilde{P}_{L}(\omega,\sigma)$.

Подставляя образ Фурье функции Лоренца (\ref{FGauss_FLorentz}) в формулу (\ref{PKR0}), получим
\begin{equation*}
P_{L}(\omega,\sigma)=\sigma^{2}\pi^{2}\sum \limits_{k=-\infty}^{\infty} e^{-2\sigma|\omega+2\pi k|}.
\end{equation*}
Раскроем модули, считая, что  $\omega \in [0,2\pi]$. Тогда
$$P_{L}(\omega,\sigma)=\sigma^{2}\pi^{2}\left(\sum \limits_{k=-\infty}^{-1} e^{2\sigma(\omega+2\pi k)}+\sum \limits_{k=0}^{\infty} e^{-2\sigma(\omega+2\pi k)}\right)=$$
$$=\sigma^{2}\pi^{2}\left(e^{2\sigma\omega}\sum \limits_{k=1}^{\infty} e^{-4\sigma\pi k}+e^{-2\sigma\omega}\sum \limits_{k=0}^{\infty} e^{-4\sigma\pi k}\right).$$
Получившиеся в последней формуле ряды легко вычисляются как суммы бесконечно убывающих геометрических прогрессий:
$$P_{L}(\omega,\sigma)=\sigma^{2}\pi^{2}\left( \frac{e^{2\sigma\omega-4\sigma\pi}}{1-e^{-4\sigma\pi}}+\frac{e^{-2\sigma\omega}}
{1-e^{-4\sigma\pi}}\right)=
\frac{\sigma^{2}\pi^{2}\ch(2\sigma(\omega-\pi))}{\sh(2\sigma\pi)}.$$
Минимальное значение $\ch(2\sigma(\omega-\pi))$ достигается при $\omega=\pi$, а максимальное -- при $\omega=0$ и $\omega=2\pi$. Следовательно
\begin{equation*}
\begin{array}{cccc}
A_{L}(\sigma)=P_{L}(\pi,\sigma)=\frac{\sigma^{2}\pi^{2}}{\sh(2\sigma\pi)}, \;
B_{L}(\sigma)=P_{L}(0,\sigma)=\frac{\sigma^{2}\pi^{2} \ch(2\sigma\pi)}{\sh(2\sigma\pi)}.
\end{array}
\end{equation*}
Теорема доказана.

\subsection{Доказательство теоремы 2}
\

Вычислим $\Phi_{L}(t,\sigma)$, пользуясь формулой (\ref{PLorentz}):

$$\Phi_{L}(t,\sigma)=\sum_{k=-\infty}^\infty \frac{\sigma^{2}}{\sigma^{2}+k^{2}} e^{-ikt}= \sigma\pi \sum \limits_{k=-\infty}^{\infty}e^{-\sigma |t+2\pi k|}.$$
Разбив, как и при доказательстве теоремы 1, последний ряд в сумму двух геометрических прогрессий, получим
\begin{equation}
\Phi_{L}(t,\sigma)=\sigma\pi \frac{\ch(\sigma(t-\pi))}{\sh(\sigma \pi)}, \, t \in [0,2\pi].
\label{PhiL}
\end{equation}
Так как $\Phi_{L}(t,\sigma)>0$ при всех значениях $t$, то на основании формулы (\ref{Int8})

$$D_{L}(t,\sigma)=\frac{1}{\Phi_{L}(t,\sigma)}=\frac{\sh(\sigma \pi)}{\sigma\pi \ch(\sigma(t-\pi))}=\sum \limits_{k=-\infty}^{\infty}d_{L,k}(\sigma) \cdot e^{-ikt}.$$
Отсюда находим
\begin{equation*}
d_{L,k}(\sigma)=\frac{\sh(\sigma\pi)}{2\sigma\pi^{2}}\int \limits_{0}^{2\pi}\frac{e^{ikt}}{\ch(\sigma(t-\pi))}dt.
\end{equation*}
Преобразуем последнюю формулу, сделав в интеграле замену переменных $t-\pi=x$
\begin{equation*}
d_{L,k}(\sigma)=\frac{\sh(\sigma\pi)\cdot e^{ik\pi}}{2\sigma\pi^{2}}\int \limits_{-\pi}^{\pi}\frac{e^{ikx}}{\ch(\sigma x)}dx=\frac{\sh(\sigma\pi)\cdot (-1)^{k}}{\sigma\pi^{2}}\int \limits_{0}^{\pi}\frac{\cos(kx)}{\ch(\sigma x)}dx.
\end{equation*}
Теорема доказана.

\subsection{Доказательство теоремы 3:\\ случай функции Лоренца}
\

Сначала найдем образ Фурье узловой функции:
\begin{equation*}
\widehat{\widetilde{\varphi}}_{L}(\omega,\sigma)=\sum_{k=-\infty}^\infty d_{L,k} \cdot
\widehat{\varphi}_{L}(\omega,\sigma)e^{-ik\omega}=
\sigma \sqrt{\frac{\pi}{2}}e^{-\sigma |\omega|}\sum_{k=
=-\infty}^{\infty} d_{L,k}(\sigma)e^{-ik\omega}.
\end{equation*}
Следовательно
\begin{equation}
\widehat{\widetilde{\varphi}}_{L}(\omega,\sigma)=
\sigma \sqrt{\frac{\pi}{2}}e^{-\sigma |\omega|}\cdot D_{L}(\omega,\sigma)=
\sigma \sqrt{\frac{\pi}{2}}e^{-\sigma |\omega|} \frac {1} {\Phi_{L}(\omega,\sigma)} \, .
\label{FUL}
\end{equation}
Вычислим $\widetilde{P}_{L}(\omega,\sigma)$ с помощью (\ref{PKR0}):
$$\widetilde{P}_{L}(\omega,\sigma)=2\pi \sum \limits_{k=-\infty}^{\infty}\left|\frac{\sigma \sqrt{\frac{\pi}{2}}e^{-\sigma |\omega+2\pi k|}}{\Phi_{L}(\omega+2\pi k,\sigma)}\right|^{2}=\sigma^{2}\pi^{2}\sum \limits_{k=-\infty}^{\infty}\frac{e^{-2\sigma |\omega+2\pi k|}}{|\Phi_{L}(\omega+2\pi k,\sigma)|^{2}} \, .$$
В силу периодичности $\Phi_{L}(\omega,\sigma)$ и с учетом (\ref{PhiL}) получим:
\begin{equation}
\widetilde{P}_{L}(\omega,\sigma)=\frac{\sigma^{2}\pi^{2}}{|\Phi_{L}(\omega,\sigma)|^{2}}\sum \limits_{k=-\infty}^{\infty}e^{-2\sigma |\omega+2\pi k|}=\frac{\sh^{2}(\sigma \pi)}{\sh(2\sigma \pi)}\cdot\frac{\ch(2\sigma(\omega-\pi))}{\ch^{2}(\sigma(\omega-\pi))}.
\label{PLFL}
\end{equation}
Для нахождения констант Рисса необходимо теперь найти максимум и минимум полученной функции на отрезке $\omega \in [0, 2\pi]$. Преобразуем последнее выражение:

$$\frac{\sh^{2}(\sigma\pi)}{\sh(2\sigma\pi)}\cdot\frac{\ch(2\sigma(\omega-\pi))}{\ch^{2}(\sigma(\omega-\pi))}=
\frac{\sh^{2}(\sigma\pi)}{\sh(2\sigma\pi)}\cdot\frac{2\ch^{2}(\sigma(\omega-\pi))-1}{\ch^{2}(\sigma(\omega-\pi))}=$$
$$=\frac{\sh^{2}(\sigma\pi)}{\sh(2\sigma\pi)}\cdot \left(2-\frac{1}{\ch^{2}(\sigma(\omega-\pi))}\right).$$
Полученная функция минимальна при $\omega=\pi$ и максимальна при $\omega=0$ или $\omega=2\pi$. Таким образом, получаем

$$
\widetilde{A}_{L}(\sigma)=\frac{\sh^{2}(\sigma\pi)}{\sh(2\sigma\pi)}, \;
\widetilde{B}_{L}(\sigma)=\frac{\sh^{2}(\sigma\pi)}{\sh(2\sigma\pi)}\left(2-\frac{1}{\ch^{2}(\sigma\pi)}\right).
$$
Осталось перейти к пределу при $\sigma\rightarrow\infty$:
$$\lim \limits_{\sigma\rightarrow\infty} \widetilde{A}_{L}(\sigma)=
\lim \limits_{\sigma\rightarrow\infty}\frac{\sh(\sigma\pi)}{2 \ch(\sigma\pi)}=
\lim \limits_{\sigma\rightarrow\infty}\frac{\th(\sigma\pi)}{2}=\frac{1}{2}.$$
$$\lim \limits_{\sigma\rightarrow\infty} \widetilde{B}_{L}(\sigma)=\lim \limits_{\sigma\rightarrow\infty}\left(\frac{\sh^{2}(\sigma\pi)}{\sh(2\sigma\pi)}\right)\cdot\lim \limits_{\sigma\rightarrow\infty}\left(2-\frac{1}{\ch^{2}(\sigma\pi)}\right)=1.$$
Теорема доказана.

\subsection{Доказательство теоремы 3:\\ случай функции Гаусса}
\

В работе \cite{Zhur} для функции $P_{G}(\omega,\sigma)$ получено соотношение
\begin{equation}
P_{G}(\omega,\sigma)= \frac { \sum \limits_{k=-\infty}^{\infty} \exp \left(- \sigma^2 (\omega+2 \pi k)^2 \right) }
{ \left[ \sum\limits_{k=-\infty}^{\infty} \exp \left( - \frac{\sigma^2}{2} (\omega+2 \pi k)^2 \right) \right] ^2} \, .
\label{PGWS_0}
\end{equation}
В этой же работе показано, что максимум $P_{G}(\omega,\sigma)$ достигается в точке $\omega=0$, вычислен предел при $\sigma \rightarrow\infty$ верхней константы Рисса
$$
\lim \limits_{\sigma\rightarrow\infty}\widetilde{B}_{G}(\sigma)=
\lim \limits_{\sigma\rightarrow\infty}\widetilde{P}_{G}(0,\sigma)=1,
$$
и получена оценка сверху для предельного значения нижней константы Рисса:
$$
\lim \limits_{\sigma\rightarrow\infty}\widetilde{P}_{G}(\pi,\sigma)=\frac {1} {2} \; \Rightarrow
\overline{\lim_{\sigma\rightarrow\infty}}\widetilde{A}_{G}(\sigma)\leq\frac{1}{2}.
$$
Для полного доказательства теоремы 3 осталось проверить, что функция $P_{G}(\omega,\sigma)$ принимает минимальное значение в точке $\omega=\pi$. Мы докажем более общее утверждение о монотонном убывании данной функции на интервале $(0,\pi)$ и монотонном возрастании на интервале $(\pi,2\pi)$.

Преобразуем знаменатель (\ref{PGWS_0}) с помощью (\ref{PGauss}) и (\ref{Teta3})
$$
\sum \limits_{k=-\infty}^{\infty} \exp{\left(-\frac{\sigma^{2}(\omega+2\pi k)^{2}}{2}\right)}=
\frac{1}{\sigma\sqrt{2\pi}}\sum \limits_{k=-\infty}^{\infty} \exp{\left(-\frac{k^{2}}{2\sigma^{2}}\right)}\cdot e^{-ik\omega}=
$$
$$
=\frac{1}{\sigma\sqrt{2\pi}} \vartheta_3 \left( \frac {\omega}{2}, \exp(-\frac{1}{2\sigma^2}) \right).
$$
Сделав аналогичное преобразование для числителя, получим равенство
\begin{equation}
P_{G}(\omega,\sigma)=\sigma \sqrt{\pi}
\frac{\vartheta_3 \left( \frac {\omega}{2}, \exp(-\frac{1}{4\sigma^2}) \right)}
     {\left[ \vartheta_3 \left( \frac {\omega}{2}, \exp(-\frac{1}{2\sigma^2}) \right) \right]^2} \, .
\label{PGWS}
\end{equation}
Введем следующие обозначения:
$$
t=\frac {\omega}{2}, \, t \in [0, \pi]; \, q=\exp \left( -\frac{1}{4\sigma^2} \right); \,
P(t)=\frac {P_{G}(\omega,\sigma)} {\sigma \sqrt{\pi}}.
$$
Тогда формула (\ref{PGWS}) примет вид
\begin{equation}
P(t)=\frac{\vartheta_3 (t, q)}
     {\left[ \vartheta_3 \left( t,q^2 \right) \right]^2} \, .
\label{PGT}
\end{equation}
Воспользуемся тождеством Ватсона \cite [c. 531] {NIST}
$$
\vartheta_3 (t, q) \vartheta_3 (w, q) = \vartheta_3 (t+w, q^2) \vartheta_3 (t-w, q^2) + \vartheta_2 (t+w, q^2) \vartheta_2 (t-w, q^2).
$$
Выбрав $w=0$, придем к равенству
$$
\vartheta_3 (t, q) \vartheta_3 (0, q) = \vartheta_3^2 (t, q^2)  + \vartheta_2^2 (t, q^2),
$$
откуда следует другая форма записи (\ref {PGT})
\begin{equation*}
P(t)=\frac{\vartheta_3 (t, q)}{\vartheta_3^2 (t, q^2)} = \frac{1}{\vartheta_3 (0, q)} \left[ 1 + \left( \frac{\vartheta_2 (t, q^2)}{\vartheta_3 (t, q^2)} \right)^2 \right].
\end{equation*}
Обозначим $p=q^2$.  Покажем, что функция $\vartheta_2 (t, p) / \vartheta_3 (t, p)$ монотонно убывает на интервале $(0,\pi)$. Согласно формуле для производной отношения двух тета-функций \cite [c. 19] {Lawden}
$$
\frac{d}{dt} \left[ \frac{\vartheta_2 (t, p)}{\vartheta_3 (t, p)} \right] = - \vartheta_4^2 (0, p)  \frac{\vartheta_1 (t, p) \vartheta_4 (t, p)}{\vartheta_3^2 (t, p)}.
$$
Из представления тета--функций в виде бесконечного произведения \cite{UW} следует, что
$\vartheta_1 (t, p)$ положительна при $t \in (0, \pi)$, а $\vartheta_3 (t, p)$ и $\vartheta_4 (t, p)$ положительны при всех $t \in \mathbf{R}$.
Следовательно, производная  $\vartheta_2 (t, p) / \vartheta_3 (t, p)$ по $t$ отрицательна при всех $t \in (0,\pi)$, а сама функция убывает. Так как $\vartheta_2 (t, p)$ положительна при $t \in (0, \frac {\pi} {2})$ и отрицательна при $t \in (\frac {\pi} {2}, \pi)$, то из формулы (\ref {PGT}) следует  монотонное убывание $P(t)$  на $(0, \frac {\pi} {2})$ и монотонное возрастание на $(\frac {\pi} {2}, \pi)$. Теорема доказана.

\subsection{Доказательство теоремы 4}
\

Так как преобразование Фурье (\ref{FT}) унитарно в $L_2(\mathbf{R})$, а образ Фурье функции $\mathrm{sinc}(\pi t)$ равен $\frac{1}{\sqrt{2\pi}}\chi_{[-\pi,\pi]}(\omega)$, где
$\chi_{[a,b]}(\omega)$ -- характеристическая функция отрезка $[a,b]$,
то доказательство теоремы удобно проводить в образах Фурье. Согласно равенству Парсеваля и формуле (\ref{FUL})
\begin{equation*}
||\widetilde{\varphi}_{L}(t,\sigma)-\mathrm{sinc}(\pi t)||^{2}_{L_{2}}=||\widehat{\widetilde{\varphi}}_{L}(\omega,\sigma)-
\frac{1}{\sqrt{2\pi}}\chi_{[-\pi,\pi]}(\omega)||^{2}_{L_{2}}=
\end{equation*}
\begin{equation}
=\int \limits_{-\infty}^{\infty}\left| \frac{\sigma \sqrt{\frac{\pi}{2}}e^{-\sigma |\omega|}}{\Phi_{L}
(\omega,\sigma)} - \frac{1}{\sqrt{2\pi}}\chi_{[-\pi,\pi]}(\omega)\right|^{2}d\omega.
\label{4start}
\end{equation}
Возведем подынтегральное выражение в квадрат. Обозначим интеграл (\ref{4start}) через $I$ и разобьем его на три слагаемых:
$I=I_{1}(\sigma)-2I_{2}(\sigma)+I_{3}(\sigma)$, где
\begin{equation*}
I_{1}(\sigma)=\frac{\sigma^{2}\pi}{2}\int \limits_{-\infty}^{\infty}\left(\frac{e^{-\sigma |\omega|}}{\Phi_{L}
(\omega,\sigma)}\right)^{2}d\omega;
\end{equation*}
\begin{equation*}
I_{2}(\sigma)=\frac{\sigma}{2}\int \limits_{-\infty}^{\infty}\frac{e^{-\sigma |\omega|}}{\Phi_{L}
(\omega,\sigma)}\cdot \chi_{[-\pi,\pi]}(\omega)d\omega;
\end{equation*}
\begin{equation*}
I_{3}(\sigma)=\frac{1}{2\pi}\int \limits_{-\infty}^{\infty}\left|\chi_{[-\pi,\pi]}(\omega)\right|^{2}d\omega=1.
\end{equation*}

Для вычисления $I_{1}(\sigma)$ разобьем $\mathbf{R}$ на отрезки $[2\pi n, 2\pi (n+1)], n \in \mathbf{Z}$:
\begin{equation*}
I_{1}(\sigma)=\frac{\sigma^{2} \pi}{2}\sum \limits_{n=-\infty}^{\infty} \int \limits_{2\pi n}^{2\pi (n+1)}\frac{ e^{-2\sigma |\omega|}}{\left(\Phi_{L}
(\omega,\sigma)\right)^{2}}d\omega.
\label{I1sum}
\end{equation*}
В каждом из интегралов сделаем замену переменных $\xi=\omega-2\pi n$. С учетом того, что знаменатель является $2\pi$-периодической функцией, получим
\begin{equation*}
I_{1}(\sigma)=\frac{\sigma^{2} \pi}{2}\sum \limits_{n=-\infty}^{\infty} \int \limits_{0}^{2\pi}\frac{e^{-2\sigma |\xi+2\pi n|}}{\left(\Phi_{L}
(\xi,\sigma)\right)^{2}}d\xi.
\label{I1sum2}
\end{equation*}
Так как функция в знаменателе $\Phi_{L}
(\xi,\sigma)$ непрерывна и строго положительна, а числитель экспоненциально убывает, то операции суммирования и интегрирования можно поменять местами. Получающуюся бесконечную сумму мы уже вычисляли в пункте 3.1 (см. формулу (\ref{PLFL})):
\begin{equation}
I_{1}(\sigma)=\frac{1}{2\pi}\frac{\sh^{2}(\sigma \pi)}{\sh(2\sigma \pi)}\int \limits_{0}^{2\pi}\frac{\ch(2\sigma(\xi-\pi))}{\ch^{2}(\sigma(\xi-\pi))}d\xi.
\label{I1sum3}
\end{equation}
Интеграл (\ref{I1sum3}) с использованием формулы для гиперболического косинуса двойного угла вычисляется аналитически. Приведем окончательный результат:
\begin{equation*}
I_{1}(\sigma)=\left(1-\frac{1}{2\sigma\pi}\cdot \th(\sigma\pi)\right)\cdot \th(\sigma\pi).
\end{equation*}
Поскольку $\th(\sigma\pi)\rightarrow 1$ при $\sigma\rightarrow \infty$, то $I_{1}(\sigma)\rightarrow 1$.

Осталось рассмотреть слагаемое $I_{2}(\sigma)$. Из-за наличия множителя $\chi_{[-\pi,\pi]}(\omega)$ пределы интегрирования становятся конечными. Кроме того, подынтегральная функция четная, поэтому
\begin{equation}
I_{2}(\sigma)=\sigma\int \limits_{0}^{\pi}\frac{e^{-\sigma \omega}}{\Phi_{L}
(\omega,\sigma)}d\omega
=\frac{\sh(\sigma\pi)}{\pi}\int \limits_{0}^{\pi}\frac{e^{-\sigma\omega}}{\ch(\sigma(\omega-\pi))}d\omega.
\label{I2hyper}
\end{equation}
Интеграл (\ref{I2hyper}) также легко вычисляется аналитически:
\begin{equation}
I_{2}(\sigma)=(1-e^{-2\sigma\pi})\cdot\left(1+\frac{1}{2\sigma\pi}\ln{\left(\frac{1+e^{-2\sigma\pi}}{2}\right)}\right).
\label{I2total}
\end{equation}
Оба сомножителя в (\ref{I2total}) стремятся к $1$ при $\sigma\rightarrow\infty$,  следовательно, и $I_{2}(\sigma)\rightarrow 1$.
Таким образом, окончательно получаем:
\begin{equation*}
\lim \limits_{\sigma \rightarrow\infty}||\widetilde{\varphi}_{L}(t,\sigma)-\mathrm{sinc}(\pi t)||^{2}_{L_{2}}=\lim \limits_{\sigma \rightarrow\infty}I_{1}(\sigma)-2\lim \limits_{\sigma \rightarrow\infty}I_{2}(\sigma)+\lim \limits_{\sigma \rightarrow\infty}I_{3}(\sigma)=0.
\end{equation*}
Теорема доказана.

\newpage
\section{Обсуждение результатов}
\

Системы целочисленных сдвигов функции Гаусса широко используются в различных разделах физики, например, в теории когерентных состояний \cite{Per}. Особенно популярной эта тематика стала после появления пакета прикладных программ ”Gaussian”, предназначенного для расчета сложных молекул. В качестве базисных функций в этом пакете используются произведения сдвигов функции Гаусса на многочлены невысоких степеней. Как показано в монографии В.Г. Мазьи, Г. Шмидта \cite{Maz},    системы сдвигов функции Гаусса могут быть применены для аппроксимации различных потенциалов, а также для решения линейных и нелинейных граничных задач математической физики.

Отдельным направлением является численная реализация рассматриваемых теоретических методов в качестве эффективных вычислительных алгоритмов. Безусловно, наличие явных формул из \cite{Maz}  в случае системы целочисленных сдвигов функций Гаусса, для которой узловая функция и коэффициенты разложений выражаются через тета--функции Якоби, имеют важное теоретическое значение. Но представляется, что основанные на подобных формулах методы не имеют вычислительных перспектив. Причина этого в том, что указанные формулы требуют деления на определённые тета--функции Якоби, которые на рассматриваемом промежутке хоть и отличны от нуля, но принимают крайне малые значения. Точнее, как показано в \cite{MS1}, при разумных значениях параметров минимумы тета--функций имеют порядки $10^{-1000}$ и меньше, что обесценивает любые попытки вычислений по явным формулам. Поэтому надо искать другие пути.

Один из вариантов вычислений для данного класса задач разработан в \cite{ZhKMS},\ \cite{SMZ}.
Он основан на применении дискретного преобразования Фурье и формул преобразования Пуассона для тета--функций. К достоинствам данного метода надо отнести то, что при его реализации, требующей вычислений с величинами порядка $10^{100}$ и больше,  удалось достигнуть точности вычислений порядка $10^{-17}$ в узлах с небольшими номерами!
Из ограничений этого метода надо отметить существенное возрастание объёмов вычислений при увеличении точности и заметные величины ошибок для приближений в узлах с номерами порядка нескольких десятков и больше.

Другой подход для вычислений был разработан в \cite{ST1}--\cite{ST6}, он является наиболее простым и прямолинейным. При этом подходе бесконечномерные системы уравнений для нахождения коэффициентов узловой функции урезаются до конечномерных систем, а затем приближаются решениями конечномерных систем. В указанных работах доказана корректная разрешимость приближающих конечномерных систем, проведены многочисленные расчёты и анализ ошибок. Этот вычислительный метод в настоящее время представляется для разложений по системам целочисленных сдвигов функций Гаусса одним из наиболее простых, но в то же время перспективных.

Для функции Лоренца системы целочисленных сдвигов изучены существенно меньше. Но, как следует из доказанных в настоящей работе теорем, предельное поведение узловых функций и связанных с ними констант Рисса для функций Гаусса и Лоренца практически одинаково. По-видимому, речь идет о какой-то общей закономерности, присущей достаточно широкому классу функций, порождающих системы целочисленных сдвигов.

Отметим, что в работах \cite{Vin1}--\cite{Vin2}  изучаются интерполяционные разложения произвольной функции по целочисленным сдвигам обобщённого распределения Коши--Лоренца вида
$$
\sum\limits_{k=-\infty}^{\infty} f_k \frac{a}{b+(x-k)^c},
$$
где $a,b,c$---некоторые положительные постоянные. Cлучай $c=2$ рассмотрен в настоящей работе.

Необходимо специально отметить, что неравенства для  тета–функций Якоби и их различных отношений играют важную роль в математических и прикладных задачах. Например, в работе \cite{Sol} неравенство о монотонности некоторого отношения тета--функций $\vartheta_2 (x, q)$ является основным моментом при решении методом диссимметризации В.\,Н. Дубинина обобщённой задачи А.\,А. Гончара о гармонической мере радиальных разрезов. Эта задача имеет длительную историю и своё развитие, в этом направлении см. \cite{Dixit} -- \cite{Schief} и библиографию в этих работах. О других неравенствах для близких спецфункций см. \cite{S1}--\cite{S6}.

Интересной и до конца не исследованной проблемой являются вопросы полноты рассматриваемых систем, это один из центральных вопросов теории когерентных состояний. Для системы когерентных состояний известно утверждение фон Неймана о полноте, в которое все верят, но которое так и не было никогда доказано. Прямое доказательство неполноты при определённых соотношениях между параметрами и обсуждение истории вопроса см. в \cite{SMZ}.

Две рассмотренные системы функций являются характерными для атомных спектров \cite [c. 23, 25] {Drob}.
Дело в том, что даже в случае дискретных спектров действие различных механизмов приводит к уширению, т.е. к образованию некоторого спектрального распределения интенсивности вблизи частоты квантового перехода в атоме или молекуле.
Функции Лоренца появляются, когда мал эффект Доплера. При сильном доплеровском уширении возникают функции Гаусса.

Сложные спектры состоят из суммы нескольких разнесенных по частотной оси функций. Требуется определить частоту и амплитуду каждой компоненты. Реальный сигнал задается, как правило, на равномерной сетке точек. Следовательно, математически задача представляет собой разложение сигнала по системе сдвигов одной функции.

Для надежного определения положения спектральной линии необходимо, чтобы на пике, связанном с этой линией, укладывалось как минимум $5-10$ дискретных отсчетов, что соответствует значениям $\sigma$ порядка $3-5$. Следовательно, анализ устойчивости с ростом параметра $\sigma$ используемых для разложения систем функций представляет большой практический интерес.

Авторами рассмотрены также некоторые другие задачи, не вошедшие в данную статью. К ним относятся вопросы оптимизации разложений по системе целочисленных сдвигов функций Гаусса для случая, когда данные представлены суммой нескольких пиков, а также приложения к теории преобразования Радона и медицинской томографии.

\newpage

В заключение сформулируем некоторые нерешённые задачи и проблемы.

Первые три задачи связаны с разложениями по сдвигам функций Гаусса.

Задача 1. Доказать, что конечномерные приближения методом работ \cite{ST1}--\cite{ST5} для нахождения узловой функции  (\ref{Int3}) сходятся к точному решению из \cite{Maz} при увеличении размерности систем.

Задача 2. Доказать, что величины коэффициентов узловой функции (\ref{Int3}) симметричны относительно нулевого и являются знакочередующимися. Некоторые результаты получены в \cite{MSU}.

Задача 3. Доказать, что экстремумы узловой функции (\ref{Int3}) достигаются при полуцелых значениях переменной.

Вычисления подтверждают истинность утверждений задач 1--3.

Задача 4. Выразить через известные специальные функции коэффициенты разложения узловой функции и саму узловую функцию для системы сдвигов функции Лоренца. Некоторые результаты получены в \cite{KMNS1}--\cite{KMNS2}.

Задача 5. Распространить полученные результаты на системы сдвигов не с целыми, а с \textsl{полуцелыми} узлами интерполяции.

Задача 6. Распространить полученные результаты на системы сдвигов с неравномерно расположенными узлами интерполяции.

Задача 7. Описать функциональные классы единственности для разложений по изученным системам.

Задача 8. Описать связь разложений по квадратичным экспонентам и когерентным состояниям с различными типами преобразования Фурье (в том числе модифицированными ДПФ \cite{S17}) и операторами свёртки, а также другими интегральными преобразованиями.

Задача 9. Получить эффективные соотношения для нахождения констант Рисса в случае двухпараметрических систем сдвигов, например, для оконного преобразования Фурье \cite{NPS} или дискретных подсистем когерентных состояний \cite{Per}.

Задача 10. Представляется важным найти применения теории операторов преобразования (см. \cite{S7}--\cite{S16}) к данному кругу задач. Наиболее интересной является задача о построении в явном виде операторов преобразования, сплетающих квадратичные экспоненты (функции Гаусса) и подобные им функции базиса когерентных состояний.

\newpage


\end{document}